\documentclass[12pt]{amsart}

\usepackage{amssymb}
\usepackage{epsf}
\usepackage{amsmath,epsf}
\usepackage[latin1]{inputenc}
\usepackage{amsfonts}

\numberwithin{equation}{section}
\newtheorem{thm}{Theorem}[section]
\newtheorem{prop}[thm]{Proposition}

\newtheorem{lem}[thm]{Lemma}
\newtheorem{cor}[thm]{Corollary}

\newtheorem{rem}[thm]{Remark}

\newcommand{\twolines}[2]{{\scriptstyle #1 \atop \scriptstyle
\vphantom{\sum}#2}}

\def\Q{{\mathbb Q}}
\def\X{{X_n}}
\def\Y{{Y_n}}

\def\S{{\mathcal S}}
\def\RS#1{{\mathbb S}_{#1}}

\def\Young#1{\vbox{\smallskip\offinterlineskip
    \halign{&\vbox{##}\kern-\Thickness\cr #1}}}

\newdimen\Squaresize \Squaresize=20pt
\newdimen\Thickness \Thickness=.1pt
\newdimen\Correction \Correction=7pt

\def\Vide#1{\hbox{
       \vbox to \Squaresize{\vss
          \hbox to \Squaresize{\hss#1 \hss}\vss}
    \hskip-\Correction}
   \kern-\Thickness}

\def\Carre#1{\hbox{\vrule width \Thickness
   \vbox to \Squaresize{\hrule height \Thickness\vss
      \hbox to \Squaresize{\hss$\scriptstyle#1$\hss}
   \vss\hrule height\Thickness}
   \unskip\vrule width \Thickness}
   \kern-\Thickness}

\title[Vanishing ideals of Lattice Diagram determinants]{Vanishing ideals
of Lattice Diagram
determinants}

\author[J.-C.~Aval and N.~Bergeron]
{J.-C.~Aval and N.~Bergeron}

\address[Jean-Christophe Aval]
{Laboratoire A2X\\ Universit\'e Bordeaux 1\\ 351 cours de la
Lib\'eration\\ 33405 Talence cedex\\ FRANCE}

\address[Nantel Bergeron]
{Department of Mathematics and Statistics\\ York University\\
   To\-ron\-to, Ontario M3J 1P3\\ CANADA}

\email[Jean-Christophe Aval]{aval@math.u-bordeaux.fr}
\email[Nantel Bergeron]{bergeron@mathstat.yorku.ca}
\urladdr[Nantel Bergeron]{http://www.math.yorku.ca/bergeron}

\date{\today}
\thanks{N. Bergeron is supported in part by NSERC}


\begin{document}

\begin{abstract}
A lattice diagram is a finite set $L=\{(p_1,q_1),\ldots ,(p_n,q_n)\}$ of
lattice cells  in
the positive quadrant. The corresponding lattice diagram determinant is
$\Delta_L(\X;\Y)=\det \|\, x_i^{p_j}y_i^{q_j}\, \|$. The
space $M_L$ is the space spanned by all partial derivatives of
$\Delta_L(\X;\Y)$. We denote by $M_L^0$ the $Y$-free component of
$M_L$. For $\mu$ a partition of $n+1$, we
denote by $\mu/ij$ the diagram obtained by removing the cell $(i,j)$
from the Ferrers diagram of $\mu$. Using homogeneous partially symmetric
polynomials, we give here a
dual description of the vanishing ideal  of the space $M_\mu^0$ and we give
the first known
description of the vanishing ideal of
$M_{\mu/ij}^0$.
 \end{abstract}

\maketitle

\section{Introduction}

A {\sl lattice diagram} is a finite set $L=\{(p_1,q_1),\ldots
,(p_n,q_n)\}$ of lattice cells  in
the positive quadrant. Following the definitions and conventions of
\cite{Berg et
al,nantel_adriano}, the coordinates $p_i\ge 0$ and $q_i\ge 0$ of a cell
$(p_i,q_i)$ indicate the row and
column  position,  respectively, of the cell in the positive quadrant.
For $\mu_1\geq \mu_2\geq \cdots \geq \mu_k>0$, we say that $\mu=(\mu_1,
\mu_2, \ldots,
\mu_k)$ is a {\sl partition} of $n$ if $n=\mu_1+\cdots +\mu_k$. We associate to
a partition $\mu$ the following lattice (Ferrers) diagram $\{(i,j)\,:\,0\le
i\le k-1,\,
0\le j\le \mu_{i+1}-1\}$ and we use the symbol
$\mu$ for both the partition and its associated Ferrers diagram.
For example, given the partition $(4,2,1)$, its Ferrers diagram is:
  $$\Young{\Carre{2,0}\cr
           \Carre{1,0}&\Carre{1,1}\cr
           \Carre{0,0}&\Carre{0,1}&\Carre{0,2}&\Carre{0,3}\cr
           }\quad.$$
This consists of the lattice cells
$\{(0,0),(1,0),(2,0),(0,1),(1,1),(0,2),(0,3)\}$.
We list the cells in {\sl lexicographic} order looking at the
second-coordinate first.
That is:
   \begin{equation}\label{lex}
    (p_1,q_1)<(p_2,q_2) \quad\iff\quad q_1<q_2\quad\hbox{or}
     \quad [q_1=q_2\hbox{ and }p_1<p_2].
   \end{equation}

Given a lattice diagram $L=\{(p_1,q_1), (p_2,q_2),\ldots , (p_n,q_n)\}$
we define the {\sl lattice diagram determinant}
  $$
  \Delta_L(\X;\Y)= \det \left\| {{x_i^{p_j}y_i^{q_j}}\over{p_j!q_j!}}\right\|_{i,j=1}^n\,,
  $$
where $\X=x_1,x_2,\ldots,x_n$ and $\Y=y_1,y_2,\ldots,y_n$.
The determinant
$\Delta_L(\X;\Y)$ is  bihomogeneous of degree $|p|=p_1+\cdots +p_n$ in $\X$
and  degree $|q|=q_1+\cdots +q_n$  in  $\Y$. To ensure that this definition
associates a unique determinant to $L$ we require
that the list of lattice cells be ordered with the lexicographic order~\ref{lex}. The factorials will
ensure that the lattice diagram determinants behave nicely under partial derivatives.

For a polynomial $P(\X;\Y)$, the vector space spanned by all the partial
derivatives of $P$ of
all orders is denoted
${\mathcal L}_\partial[P]$. A permutation $\sigma\in \S_n$ acts diagonally on
a polynomial $P(\X;\Y)$ as follows:
$\sigma  P(\X;\Y)\,=\, P(x_{\sigma_1},x_{\sigma_2},\ldots
,x_{\sigma_n};y_{\sigma_1},y_{\sigma_2},\ldots ,y_{\sigma_n})$.
Under this action, $\Delta_L(\X;\Y)$ is clearly an alternant. Moreover, partial derivatives commute
with the action, hence it follows that for any lattice diagram $L$ with $n$ cells, the vector space
$M_L = {\mathcal L}_\partial[\Delta_L(\X;\Y)]$
is an $\S_n$-module. Since $\Delta_L(\X;\Y)$ is  bihomogeneous, this module
affords a natural
bigrading. Denoting by ${\mathcal H}_{r,s}[M_L]$ the subspace consisting of the
bihomogeneous elements of degree  $r$ in $\X$ and degree $s$ in $\Y$, we have
the direct sum
decomposition
  $$
  M_L = \bigoplus_{r=0}^{|p|}  \bigoplus_{s=0}^{|q|} {\mathcal H}_{r,s}[M_L].
  $$
In general, not much is known on the $\S_n$-modules $M_L$. In the case where $L$ is a
Ferrers diagram
$\mu$ with $n$ cells,
 the $n!$ conjecture of Garsia and Haiman \cite{GH} stated that $M_\mu$
has dimension
$n!$, that it is 
$\S_n$-isomorphic to the left regular representations and its graded character is a renormalization
of the Macdonald polynomial \cite{Mac} indexed by $\mu$. Haiman in 
\cite{haiman} has recently completed a proof of this using an algebraic geometry approach.
It develops that a very natural and combinatorial recursive approach to the result above involves
diagrams obtained by removing a single cell from a partition diagram. In
\cite{Berg et al}, we have
investigated this case and formulated
various new conjectures.

To pursue the investigation of the spaces $M_L$, we are interested in an
explicit
description of the vanishing ideal $I_L$ of differential operators on the
spaces $M_L$, which is defined in the following way:
$$I_L=\{P(\X;\Y)\in \Q[\X;\Y]:\ P(\partial \X;\partial \Y)\Delta_L(\X;\Y)=0\},$$
where for a polynomial $P(\X;\Y)$ we denote by $P(\partial \X;\partial \Y)$ the differential
operator obtained from $P$,
substituting every variable $x_i$ by the operator
$\partial\,\,\,\,\over\partial x_i$ and every variable $y_j$ by the operator
$\partial\,\,\,\,\over\partial y_j$.
A step in this program is to describe the vanishing ideal of a special
subspace of
$M_L$. For a lattice
diagram
$L$, we let
  $$
  M_L^0 = \bigoplus_{r=0}^{|p|}  {\mathcal H}_{r,0}[M_L].
  $$
That is the $Y$-free component of $M_L$.
In this paper, we study the vanishing ideal $I_L^0$ of differential
operator on the spaces $M_L^0$. 
Our result gives a good set of generators for the ideal $I_L^0$ in the
case where $L$ is a
partition $\mu$ or a partition with a hole $\mu/ij$. The case $L=\mu$ was
studied extensively in
\cite{nantel_adriano,CP,GP,Tani}. Our description is dual to 
Tanisaki's~\cite{nantel_adriano,Tani}.
For
$L=\mu/ij$, there was no previously known description of
$I_{\mu/ij}^0$.

Our analysis is based on a careful study of the effect of partially symmetric
differential
operators on lattice diagram determinants. We do this in Section~3. For
completeness, we also give a
description for Schur symmetric differential operators. In Section~4 we
recall the known results for $I_\mu^0$
and give the dual description using  homogeneous partially symmetric
polynomials. In Section~5 we describe
$I_{\mu/ij}^0$.

\section{Some preliminary results}

Before we start, let us collect some useful facts. 
We need a few definitions. 
For an $n$-cell lattice diagram $L$, a {\sl tableau} of shape $L$ is an injective map $T:L\to
\{1,2,\ldots,n\}$. We can think of $T$ as a way to list the cells of $L$. If $T(r,c)=m$, we say that
$h_T(m)=r$ is the {\sl height\/} of $m$ in $T$. We say that $T$ is {\sl column increasing} if
$T(r,c_1)<T(r,c_2)$ whenever
 $c_1<c_2$ (when this has a meaning). Let ${\mathcal T}_L$ be the set of all tableaux of shape $L$ and
let ${\mathcal CT}_L$ be the set of all column increasing tableaux of shape $L$. Finally, for any
tableau $T$ of shape $L$, we let 
 $$\Gamma(T)=\big\{ \{T(r,k) | (r,k)\in L\}\quad \big|\quad 0\le k\le\max_{(r,c)\in L} c\big\}$$
be the {\sl column sets} of $T$.
 
We now expand the determinant
  $$\Delta_L(\X,\Y) = \Big(\prod_{(r,c)\in L} {1\over r!\,c!}\Big)\sum_{T\in{\mathcal T}_L} 
    \pm \, m_T(\Y) \widetilde{m}_T(\X)$$ 
where 
  $$ m_T(\Y)=\prod_{(r,c)\in L}y^c_{T(r,c)},\quad\qquad 
     \widetilde{m}_T(\Y)=\prod_{(r,c)\in L}x^r_{T(r,c)},$$ 
and the sign is the sign of the permutation that reorders the cells of $L$ given by $T$ back in the
lexicographic order~\ref{lex}.  If we now collect the terms of $\Delta_L(\X,\Y)$ with the same
monomials in $\Y$, we easily see that
  \begin{equation} \label{introexpand}
  \Delta_L(\X,\Y) = \Big(\prod_{(r,c)\in L} {1\over r!\,c!}\Big)\sum_{T\in{\mathcal CT}_L} 
    m_T(\Y) \Delta_T(\X)
  \end{equation}
where 
  \begin{equation} \label{garnirD}
   \Delta_T(\X)=\pm\prod_{{\mathcal C}\in\,\Gamma(T)} \det\left\|
     x_m^{h_T(\ell)}
     \right\|_{m,\ell\in\,{\mathcal C}}.
  \end{equation}

For $P(\X)\in\Q[\X]$, we have
 \begin{equation}\label{introchar}  
  P(\partial\X)\Delta_L(\X;\Y)=0 \quad\iff\quad P(\X)\in I^0_L.
 \end{equation}
On one hand, if $P(\partial\X)\Delta_L(\X;\Y)=0$ then clearly
$P(\partial\X) Q(\X;\Y)=0$ for all $Q(\X;\Y)\in M_L$. In particular $P(\partial\X) Q(\X)=0$ for all
$Q(\X)\in M^0_L\subseteq M_L$ and  $P(\X)\in I^0_L$. On the other hand, consider the
expansion~\ref{introexpand} of $\Delta_L(\X;\Y)$. For any column strict tableau $T$ of shape
$L$ we have
$m_T(\partial
\Y)
\Delta_L(\X;\Y)= c\,\Delta_T(\X)$ for a non-zero constant $c$, and $\Delta_T(\X)\in
M_L^0$. Hence if $P(\X)\in I^0_L$, we have that
$P(\partial\X)\Delta_T(\X)=0$ for all $T$ in the expansion~\ref{introexpand} and thus
$P(\partial\X)\Delta_L(\X;\Y)=0$.

\section{Symmetric operators}

For the sake of simplicity, we limit our descriptions to $X$-operators; the
$Y$-operators are
similar. Recall that
\begin{itemize}
\item  $\displaystyle p_k(\X)=\sum_{i=1}^n x_i^k$,
\item  $\displaystyle e_k(\X)=\sum_{1\le i_1 < i_2 <\cdots<i_k\le n} x_{i_1}
x_{i_2}
    \cdots x_{i_k}$\quad and
\item  $\displaystyle h_k(\X)=\sum_{1\le i_1 \le i_2 \le\cdots\le i_k\le n}
x_{i_1} x_{i_2}
    \cdots x_{i_k}$
\end{itemize}
are the power sum, elementary and homogeneous symmetric polynomials,
respectively. In
\cite{aval2,Berg et al} we have described the effect of the above symmetric
differential operators
on lattice diagram determinants. We now recall the results stated
in \cite{aval2}. 

In this section we will assume that a lattice diagram is a lists
of cells 
  $$L=[(p_1,q_1), (p_2,q_2),\ldots , (p_n,q_n)]$$
 ordered by the order~\ref{lex}.
We allow $L$ to have repeated entries or negative coordinates, but in that case we set the
determinant 
$\Delta_L(\X;\Y)=0$. We define the function
 \begin{equation}\label{coef}\epsilon(L)= \left\{ 
     \begin{array}{rll}
        1&& \hbox{if $L$ has  distinct cells in the positive quadrant,}\\
      &&\\
        0&& \hbox{otherwise.}
      \end{array}
     \right.
 \end{equation}
 Let
$L$ be a lattice diagram as above with $n$ distinct cells in the positive quadrant. For
any integer $k\ge 1$ we have:

\begin{prop}[Proposition I.1 \cite{Berg et al}]
   $$p_k(\partial \X)\Delta_L(\X;\Y)=\sum_{i=1}^n\,\,\pm\Delta_{p_k(i,L)}(\X;\Y)$$
 where ${p_k(i,L)}$ is the lattice diagram obtained from $L$ by replacing
 the $i$-th biexponent $(p_i,q_i)$ with $(p_i-k,q_i)$, and the sign of $\Delta_{p_k(i,L)}(\X;\Y)$  is
the sign of the
 permutation that reorders the resulting biexponents in the order~\ref{lex}.
\end{prop}

\begin{prop}[Proposition 2 \cite{aval2}]\label{elem}
   $$e_k(\partial \X)\Delta_L(\X;\Y)=\sum_{1\le i_1<i_2<\cdots<i_k\le
    n}\Delta_{e_k(i_1,\ldots,i_k;L)}(\X;\Y)$$ 
where
$e_k(i_1,\ldots,i_k;L)$ is the lattice diagram obtained from $L$ by
replacing the biexponents
$(p_{i_1},q_{i_1}),\ldots,(p_{i_k},q_{i_k})$ with
$(p_{i_1}-1,q_{i_1}),\ldots,(p_{i_k}-1,q_{i_k})$.
\end{prop}

For a lattice diagram $L$ with $n$ distinct cells in the positive quadrant, we denote by $\overline
L$ its complement in the positive quadrant
(it is an infinite subset). Again we list $\overline L=[(\overline
p_1,\overline q_1), (\overline
p_2,\overline q_2),\ldots\,]$ using the lexicographic order~\ref{lex}.

\begin{prop}[Proposition 3 \cite{aval2}]\label{homo}
  $$h_k(\partial \X)\Delta_L(\X;\Y)=\sum_{1\le i_1<i_2<\cdots<i_k}
      \epsilon\big(e_k(i_1,\ldots,i_k,\overline
      L)\big)\Delta_{\overline{e_k(i_1,\ldots,i_k,\overline L)}}(\X;\Y)$$
 where $e_k(i_1,\ldots,i_k,\overline L)$ is the lattice diagram obtained from $\overline L$ by
replacing the biexponents
$(\overline p_{i_1},\overline q_{i_1}),\ldots,(\overline
p_{i_k},\overline q_{i_k})$ with
$(\overline p_{i_1}+1,\overline q_{i_1}),\ldots,(\overline
p_{i_k}+1,\overline q_{i_k})$. Its complement is $\overline{e_k(i_1,\ldots,i_k,\overline L)}$ and
the function $\epsilon$ is defined in \ref{coef}.
\end{prop}

We remark that even if the sum in Proposition~\ref{homo} is infinite, only
a finite number of
coefficients $\epsilon\big(e_k(i_1,\ldots,i_k,\overline L)\big)$ differs from zero.
We will give at the end of this section an expression for $h_k(\partial \X)\Delta_L(\X;\Y)$
that does not depend on the complement of $L$. 
 For this, it is natural to ask  what is the effect of a Schur differential
operator on a lattice
diagram determinant. For completeness, we insert here such a description
that unifies
Proposition~\ref{elem} and~\ref{homo}. We only sketch a proof of our
description since this is not
used in the remaining sections and the technique of proofs are well known.

 Following \cite{Mac}, recall that for a partition
$\lambda=(\lambda_1,\lambda_2,\ldots,\lambda_k)$ the conjugate (transpose)
partition is denoted by
$\lambda'=(\lambda'_1,\lambda'_2,\ldots,\lambda'_\ell)$. With this in mind,
the Schur polynomial
indexed by
$\lambda$ is
 \begin{equation}\label{sexp}
   S_\lambda(\X)=\det{\|\,e_{\lambda'_i+j-i}(\X)\,\|} = \sum_{\sigma\in \S_\ell } sgn(\sigma)
   e_{\sigma(\lambda'+\delta_\ell)-\delta_\ell},
 \end{equation}
with the understanding that $e_0(\X)=1$ and $e_k(\X)=0$ if
$k<0$. The Schur
polynomials also have a description in terms of column-strict Young tableaux.
Given $\lambda$ a
partition of $n$, a tableau of shape $\lambda$ is a map
$T\colon\lambda\to\{1,2,\ldots,n\}$. We say that $T$ is a column-strict Young tableau if it
is weakly increasing
along the rows and strictly increasing along columns of $\lambda$. That is
$T(i,j)\le T(i,j+1)$ and
$T(i,j)<T(i+1,j)$. We denote by ${\mathcal T}_\lambda$ the set of all
column-strict Young tableaux of shape $\lambda$. For any tableau $T$, we define
$\X^T=\prod_{i=1}^n x_i^{|T^{-1}(i)|}$. As seen in~\cite{Mac}, we have
  $$S_\lambda(\X)=\sum_{T\in{\mathcal T}_\lambda} \X^T.$$

Let  $L$ be a lattice diagram with $n$ cells ordered by the
order~\ref{lex}. For any partition
$\lambda$ of an integer $k\ge 1$ we have

\begin{prop}\label{schur}
$$S_\lambda(\partial \X)\Delta_L(\X;\Y)=\sum_{T\in{\mathcal T}_\lambda}
  \epsilon'(T,L)\Delta_{\partial T(L)}(\X;\Y)$$
 where
$\partial T(L)$ is the lattice diagram obtained from $L$ by replacing the
biexponents
$(p_i,q_i)$ with
$(p_i-|T^{-1}(i)|,q_i)$ for $1\le i\le n$. The coefficient
  $$\epsilon'(T,L)=
    \epsilon\big(\partial T(L)\big)\ \cdots\ 
    \epsilon\big(\partial T_{\ell-1}\partial T_\ell(L)\big)\ \epsilon\big(\partial T_\ell(L)\big),
  $$
where $T_1,T_2,\ldots,T_\ell$ are the $\ell$ columns of $T$  and
 $\epsilon$ is the function defined in \ref{coef}. 
\end{prop}

There are many standard ways to prove this statement. For example one can
iteratively use the
Proposition~\ref{elem} with the expansion \ref{sexp}.
Then, Proposition~\ref{schur} follows by a suitable canceling involution.
See \cite{RS}
for a similar involution. A complete proof of this Proposition will appear elsewhere \cite{Sproof}.

\begin{rem}\rm Given a lattice diagram $L$ and a column strict tableau
$T\in\mathcal{T}_\lambda$, we have that $\epsilon'(T,L)=1$ exactly when we can {\sl move} the cells
of $L$ by one, reading $T$ column by column, from right to left, without having any cells
colliding.
\end{rem}

\begin{cor} For $h_k(\X)=s_{(k)}(\X)$ we have 
 $$h_k(\partial \X)\Delta_L(\X;\Y)=\sum_{1\le i_1\le i_2\le\cdots\le i_k\le
n}\epsilon'((i_1,\ldots,i_k),L)\Delta_{\partial_{i_1}\cdots\partial_{i_k}(L)}(\X;\Y)$$
\end{cor}

This is equivalent to the description in Proposition \ref{homo}. The only way that\break
$\epsilon'((i_1,\ldots,i_k),L)\ne 0$ is if the cells $i_1,\ldots,i_k$ that move down are moved
into holes. This can be described as distinct holes moving up.

\bigskip
For the following sections we now describe a necessary condition that tests
if a  partially symmetric operator
belongs to the vanishing ideal of a lattice diagram determinant.
For $k\le n$, fix
   $$S=\{x_{i_1},x_{i_2},\ldots,x_{i_k}\}\subseteq \X,$$
and set $S^Y=\{y_{i_1},y_{i_2},\ldots,y_{i_k}\}$.
Given $L$ a lattice diagram with $n$ cells, we can expand the determinant
$\Delta_L(\X;\Y)$ in term of the
rows ${i_1},{i_2},\ldots,{i_k}$ and we get
   \begin{equation}\label{Lag}
   \Delta_L(\X;\Y)=\sum_{\twolines{D\,\subset\, L,}{\null |D|\,=\,k}} \pm
\Delta_D(S;S^Y)
   \Delta_{L-D}(\X-S,\Y-S^Y).
   \end{equation}
Clearly, if a symmetric operator in the variables $S$ annihilates all
$\Delta_D(S;S^Y)$ for $D\subseteq L$,
then it  annihilates $\Delta_L(\X;\Y)$. For example, let
$L=\{(1,0),(0,1),(1,1),(0,2)\}$ and choose $S\subseteq \X$
such that $|S|=k=2$. We have that
$h_3(S)\in I^0_L$, since for all subsets $D\subset L$ such that $|D|=2$ the
Proposition~\ref{homo} gives us
that $h_3(\partial S)\Delta_D(S;S^Y)=0$. We can visualize this as follows.
We represent all subsets $D$ putting
two $\bullet$ in the all possible ways inside $L$.
\newdimen\Squaresize \Squaresize=12pt
  $$\Young{\Carre{\bullet}&\Carre{\bullet}\cr
           &\Carre{}&\Carre{}\cr
           }\quad,\quad
    \Young{\Carre{\bullet}&\Carre{}\cr
           &\Carre{\bullet}&\Carre{}\cr
           }\quad,\quad
    \Young{\Carre{\bullet}&\Carre{}\cr
           &\Carre{}&\Carre{\bullet}\cr
           }\quad,\quad
    \Young{\Carre{}&\Carre{\bullet}\cr
           &\Carre{\bullet}&\Carre{}\cr
           }\quad,\quad
    \Young{\Carre{}&\Carre{\bullet}\cr
           &\Carre{}&\Carre{\bullet}\cr
           }\quad,\quad
    \Young{\Carre{}&\Carre{}\cr
           &\Carre{\bullet}&\Carre{\bullet}\cr
           }
$$
The maximal number of distinct cells in $\overline D$ that can go up
without overlapping another in all
pictures is two. For example, in the first picture, there are two cells
below the two $\bullet$, if we try to
move up three cells of $\overline D$, necessarily at least two will overlap and
hence
$h_3(\partial S)\Delta_D(S;S^Y)=0$. Using this technique we can easily check
that $h_r(S)\in I_L^0$
if \  $r\ge 2$ for $|S|=1$,\quad $r\ge 3$ for $|S|=2$,\quad $r\ge 3$ for
$|S|=3$ and $r\ge 2$ for $|S|=4$.

\section{The ideal $I_\mu^0$}

In this section we recall the results \cite{nantel_adriano,Tani} for the
ideal $I_\mu^0$ and give a dual
description of the ideal in terms of homogeneous partially symmetric polynomials.
This  gives us a better understanding for the case $I_{\mu/ij}^0$ in Section~5.

Let $\mu$ be a fixed partition
of
$n$. The homogeneous partially symmetric polynomials are the polynomials
$h_r(S)$ for
$S=\{x_{i_1},x_{i_2},\ldots,x_{i_k}\}\subseteq \X$.  Let
$\mu'=(\mu'_1,\mu'_2,\ldots,\mu'_m)$ denote the conjugate partition of
$\mu$. For $1\le k\le
n$, we define
  $$\delta_k(\mu)=\mu'_1+\mu'_2+\cdots+\mu'_k
$$
with the convention that $\mu'_j=0$ if $j>m$.

\begin{prop}\label{dtani} For $\mu$ a partition of $n$
  $$
    I_\mu^0 = \big\langle h_r(S) : S\subseteq \X,\ |S|=k,\ r>\delta_k(\mu)-k
\big\rangle .
  $$
\end{prop}

 Tanisaki \cite{Tani}, with a simplified proof in \cite{nantel_adriano},
shows that
  $$\begin{array}{rcl}
   I_\mu^0 &=&\big\langle e_r(\overline S) : \overline S\subseteq \X,\
|\overline S|=\overline k,\
r>\overline k-\big(n-\delta_{n-\overline k}(\mu)\big)
\big\rangle
\\
             &=&\big\langle e_r(\overline S) : \overline S\subseteq \X,\
|\overline S|=n-k,\
                r>\delta_k(\mu)-k \big\rangle .
    \end{array}
$$
In light of the following lemma, the Proposition~\ref{dtani} is simply the
dual description of Tanisaki's description of $I_\mu^0$.
More precisely, using the ideal
  $$I_{1^n}^0 =\langle e_r(\X) : r>0 \rangle = \langle h_r(\X) : r>0 \rangle
\subseteq I_\mu^0,
$$
we have:

\begin{lem}\label{htoe} For $S\subseteq \X$ let ${\overline S}=\X-S$, then
 $$h_r(S)\equiv (-1)^r e_r(\,{\overline S}\,) \quad\mod\ I_{1^n}^0  \ . $$
\end{lem}

\proof We know that $e_r(\,{\overline S}\,)$ and $h_r(S)$ are the
coefficient of $t^r$ in
  $$E_{\overline S}(t)=\prod_{i\in{\overline S}} (1+tx_i)\qquad\hbox{and}\qquad
    H_{S}(t)=\prod_{i\in S} {1\over 1-tx_i}
   $$
respectively. Since $S\cup {\overline S}=\X$ is a disjoint union, we have
that $E_{\overline
S}(t)/H_{S}(-t)=E_\X(t)\equiv 1\mod I_{1^n}^0$. Hence $E_{\overline
S}(t)\equiv H_{S}(-t)\mod I_{1^n}^0$ and the result follows.
\endproof

\begin{rem}\rm The reader should note that the argument of
\cite{nantel_adriano,Tani} could be further
simplified using directly the generators of Proposition~\ref{dtani}. As
described at the end of Section 3,
one can easily show that $h_r(\partial S)\Delta_\mu(\X;\Y)=0$ whenever
$|S|=k$ and $r>\delta_k(\mu)-k$.
This is a simple use of the pigeon hole principle. Then the reduction
algorithm in the proof of Proposition
4.2~\cite{nantel_adriano} is better suited to the homogeneous
functions $h_k(S)$.
\end{rem}

\section{The ideal $I_{\mu/ij}^0$}

We describe the vanishing ideal $I_{\mu/ij}^0$ of the space $M^0_{\mu/ij}$
for $\mu$ a partition of $n+1$ and
$(i,j)$ a cell of $\mu$.
In the previous section we have seen that $I_{\mu}^0$ has two dual
descriptions in terms of 
elementary partially symmetric function and in terms of homogeneous partially symmetric
functions. The key idea
was that they are both equivalent modulo the ideal $I_{1^n}^0$ of fully
symmetric functions that is
contained in every $I_{\mu}^0$. To study $I_{\mu/ij}^0$ we do not have
$I_{1^n}^0\subseteq I_{\mu/ij}^0$.
For example, if $(i,j)\in\mu$ is not at the top of a column of $\mu$ then
$h_1(\partial \X)\Delta_{\mu/ij}(\X;\Y)=\Delta_{\mu/i+1,j}(\X;\Y)\not = 0$. 
On the other hand, as we
will see in Lemma~\ref{lem4.2},
$h_r(\X)\in I_{\mu/ij}^0$
for all other $r>1$.
To describe  $I_{\mu/ij}^0$ we need to use both families of generators.
Let us introduce some notation.

\centerline{
\epsffile{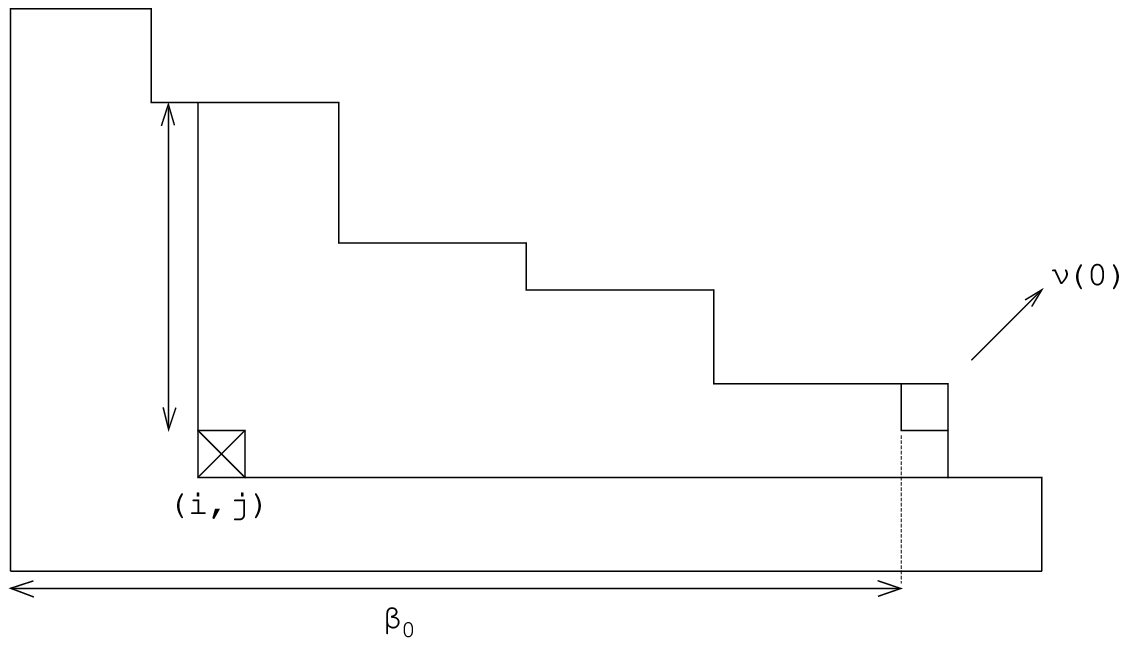}}

\centerline{\bf Figure 1}

\medskip
Let $\ell$ be the number of cells above the cell $(i,j)$ in $\mu$. In the
area north east of the cell $(i,j)$,
let $(\alpha_0,\beta_0)$ be the coordinate of the rightmost corner cell
of $\mu$. We let $\nu(0)=\mu/\alpha_0\beta_0$ be the partition of $n$
obtained from
$\mu$ by removing the cell $(\alpha_0,\beta_0)$.

If $(i,j)$ is at the top of a column of $\mu$, then $\alpha_0=i$ and
  $$h_1^{\beta_0-j}(\partial \Y) \Delta_{\mu/ij}(\X;\Y)\, =\,
   \Delta_{\mu/\alpha_0\beta_0}(\X;\Y)\,=\,\Delta_{\nu(0)}(\X;\Y).$$
If $P(\partial\X)\Delta_{\mu/ij}(\X;\Y)=0$, then clearly
$P(\partial\X)\Delta_{\nu(0)}(\X;\Y)=0$ and  $I^0_{\mu/ij}\subseteq I^0_{\nu(0)}$.
To see the other inclusion, we first compare the expansions~\ref{introexpand} for
$\Delta_{\mu/ij}(\X;\Y)$ and $\Delta_{\nu(0)}(\X;\Y)$. We have 
 $$\Delta_{\mu/ij}(\X;\Y)=\sum_{T\in{\mathcal CT}_{\mu/ij}} m_T(\Y) \Delta_T(\X)$$
 $$\Delta_{\nu(0)}(\X;\Y)=\sum_{T'\in{\mathcal CT}_{\nu(0)}} m_{T'}(\Y) \Delta_{T'}(\X)$$
There is an obvious bijection ${\mathcal CT}_{\mu/ij}\to {\mathcal CT}_{\nu(0)}$ that preserves the
column set\break $\Gamma(T)=\Gamma(T')$. Under this pairing $T\leftrightarrow T'$ we have that
$\Delta_T(\X) =
\pm
\Delta_{T'}(\X)$. If
$P(\partial\X)\Delta_{\nu(0)}(\X;\Y)=0$ we must have that
$0=P(\partial\X)\Delta_{T'}(\X)=P(\partial\X)\Delta_T(\X)$ for all $T\leftrightarrow T'$. It follows
that $P(\partial\X)\Delta_{\mu/ij}(\X;\Y)=0$ and this shows the inclusion $I^0_{\nu(0)}\subseteq
I^0_{\mu/ij}$. Hence if $(i,j)$ is at the top of a column of $\mu$, we have
$I^0_{\mu/ij}=I^0_{\nu(0)}$, a case covered in the previous section.
For the remaining of this section we can thus assume that $(i,j)$ is not at
the top of a column of $\mu$.

For $\mu$ a partition of $n+1$ let
  $$ J_\mu^0=\big\langle h_r(S) : S\subseteq \X,\ |S|=k,\ r>\delta_k(\mu)-k,\ \hbox{ for }
     1\le k\le n \big\rangle.
$$
In particular $h_1(\X)\not\in J_\mu$, but we clearly have
  $$ J_\mu^0\subseteq\big\langle h_r(S) : S\subseteq \X,\ |S|=k,\ r>\delta_k(\mu)-k,\ \hbox{ for }
     1\le k\le n+1 \big\rangle.
$$
Thus from Proposition~\ref{dtani} it follows that $J_\mu^0\subseteq I_\mu^0$.
We have seen above that the use of homogeneous partially symmetric functions
may not be enough to describe
$I_{\mu/ij}^0$. In our main theorem below, we also need elementary
partially symmetric functions. Let
  $$\begin{array}{rcl}
    \widehat J_{\mu/ij}^0 &=&\big\langle h_1(\X) e_r(\overline S) :
\overline S\subseteq \X,\ |\overline S|=n-k,\
       j<k\le\beta_0,\ r=\delta_k(\mu)-k  \big\rangle\\
        && + \quad\big\langle e_r(\overline S) : \overline S\subseteq \X,\
|\overline S|=n-k,\
       \beta_0<k<\mu_1,\ r=\delta_k(\mu)-k  \big\rangle.\\
    \end{array}$$

\begin{thm}\label{Tmuij} Using the notation above for
$\mu=(\mu_1,\mu_2,\ldots)$ a partition of $n+1$ and
$(i,j)$ a cell of
$\mu$,
   \begin{equation}\label{Imuij}
   I_{\mu/ij}^0 =  J_\mu^0\ + \ \big\langle h_1^{\ell+1}(\X) \big\rangle
+\widehat J_{\mu/ij}^0 .
   \end{equation}
\end{thm}

The remaining of this section is dedicated to a proof of this theorem.
 Let $\widetilde{I}_{\mu/ij}^0$ denote the right hand side of
Equation~\ref{Imuij}. We first show the
inclusion $\widetilde{I}_{\mu/ij}^0\subseteq I_{\mu/ij}^0$,  
showing that each of the components
of $\widetilde{I}_{\mu/ij}^0$ belongs to $I_{\mu/ij}^0$.

\begin{lem}\label{lem4.2} For any $(i,j)\in\mu$, we have $J_\mu^0\subseteq I_{\mu/ij}^0$.
\end{lem}

\proof
It is clear from the definition that $J_\mu^0\subseteq I_\mu^0$. Let us
use the Equation~\ref{Lag}
with $\mu$ and $x_{n+1}$:
    $$\Delta_\mu(\X,x_{n+1};\Y,y_{n+1})=\sum_{(i,j)\in\mu}\  \pm\, x_{n+1}^i
y_{n+1}^j\, \Delta_{\mu/ij}(\X;\Y).$$
For any $P(\X)\in J_\mu^0\subseteq I_\mu^0$ we have
    $$ 0=P(\partial \X) \Delta_\mu(\X,x_{n+1};\Y,y_{n+1})=\sum_{(i,j)\in\mu}\
\pm\, x_{n+1}^i y_{n+1}^j\,
P(\partial \X)\Delta_{\mu/ij}(\X;\Y)$$
which shows that $P(\partial \X)\Delta_{\mu/ij}(\X;\Y)=0$ for all $(i,j)\in\mu$.
\endproof

To show that $h_1^{\ell+1}(\X)\in I_{\mu/ij}^0$ we use the
Proposition~\ref{homo}. In $\overline {\mu/ij}$
the only cell that can move up, without overlapping another, is $(i,j)$. It
can move up by $\ell$, but not more.
Moving the cell by $\ell+1$ would be outside of $\mu$ hence overlapping another cell
of $\overline {\mu/ij}$. Thus
  $$h_1^{\ell+1}(\partial \X)\Delta_{\mu/ij}(\X;\Y)=0.
$$

Now, for $\beta_0<k<\mu_1$ we consider $e_r(\overline S)$ where $\overline
S\subseteq \X$ of cardinality
$|\overline S|=n-k$ and
$r=\delta_k(\mu)-k$. Again we shall show that $e_r(\overline S)\in I_{\mu/ij}^0$. We
remark that if the condition
$\beta_0<k<\mu_1$ is non-empty, then  $(i,j)$ is not on the first row of
$\mu$. We expand
 $\Delta_{\mu/ij}(\X;\Y)$ as in~\ref{Lag}
  $$
   \Delta_{\mu/ij}(\X;\Y)=\sum_{\twolines{D\,\subset\, \mu/ij,}{\null
|D|\,=\,n-k}} \pm \Delta_D(\overline S;\overline S^Y)
   \Delta_{\mu/ij-D}(\X-\overline S,\Y-\overline S^Y).
$$
By inspection of Figure~2 below, we note that there are $n-k$ cells in
total in the two shaded areas of
$\mu/ij$. In the darker grey area there are $\delta_k(\mu)-k-1$ cells of
$\mu/ij$. According to the
Proposition~\ref{elem}, to have
$e_r(\partial \overline S)\Delta_D(\overline S;\overline S^Y)\not = 0$, we must be able to move
$r=\delta_k(\mu)-k$ distinct cells of $D$
down. 
Among all the diagrams $D\,\subset\, \mu/ij$ such that $|D|\,=\,n-k$, the ones that
maximize the number of distinct cells that  can go down are obtained as follows.
We must first choose all the cells of $\mu/ij$ that are in rows $2,3,\ldots,\mu'_1$, and
there are $n-\mu_1$ such cells. We must then choose $\mu_1-k$ cells in the first row. Each cell
in the first row prevents the cells above it to move down. We must minimize these obstructions and
choose cells $(0,j_1), (0,j_2),\ldots, (0,j_{\mu_1-k})$ in columns $j_s$ such that the $\mu'_{j_s}$
are the smallest, and possibly the column $j$ (corresponding to the hole $(i,j)$). Choosing the
column $j$ would prevent $i\ge\mu_{k+1}'$ cells to move down; hence up to a permutation of the
columns, we may choose $j_s=k+s$. That is the diagram $D$ depicted by the two shaded areas in 
Figure~2. The maximal number of distinct cells that can go
down for that $D$ is
$\delta_k(\mu)-k-1$. All other $D\,\subset\, \mu/ij$ such that $|D|\,=\,n-k$ will have no more than
$\delta_k(\mu)-k-1$ distinct cells that can go down,   hence
$e_r(\partial \overline S)\Delta_{\mu/ij}(\X;\Y)=0$ and
$e_r(\overline S)\in I_{\mu/ij}^0$.

\centerline{
\epsffile{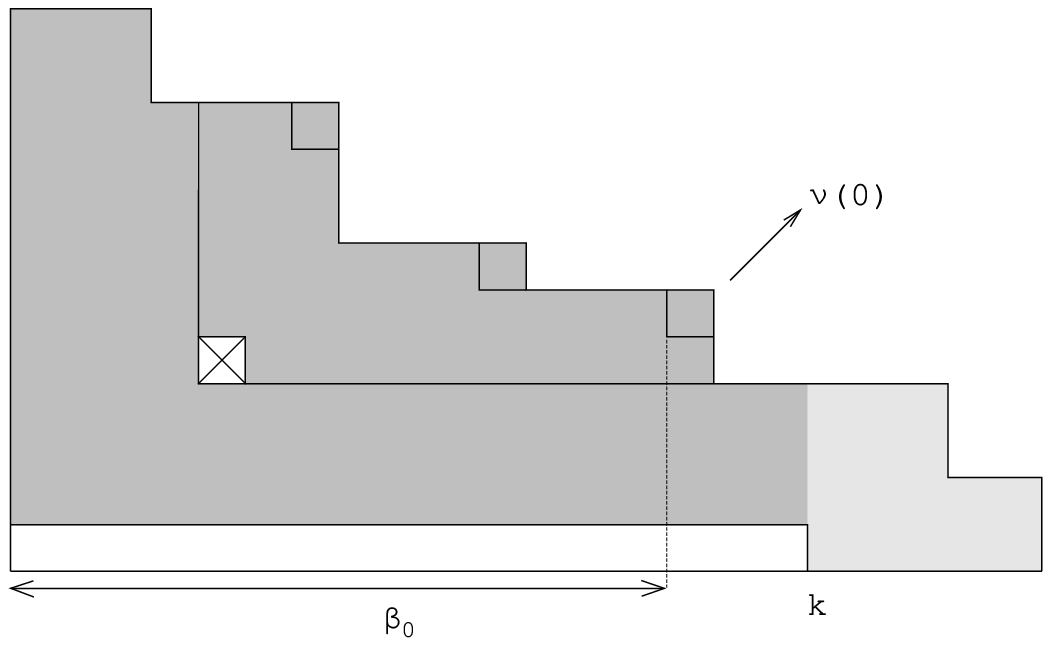}}

\centerline{\bf Figure 2}
\medskip

\noindent
Finally, for $j<k\le\beta_0$ and $e_r(\overline S)$ as above, the only remaining
problem with our argument
is when the cell $(i,j)$ is not in the darker grey area. This happens only if
$i=0$. But recall that
$h_1(\partial \X)\Delta_{\mu/ij}(\X;\Y)=\Delta_{\mu/i+1,j}(\X;\Y)$ and hence
  $$e_r(\partial \overline S)h_1(\partial
\X)\Delta_{\mu/ij}(\X;\Y)=e_r(\partial \overline
   S)\Delta_{\mu/i+1,j}(\X;\Y)=0.$$
Thus $h_1(\X)e_r(\overline S)\in I_{\mu/ij}^0$ which concludes the proof of
the inclusion
$\widetilde{I}_{\mu/ij}^0\subseteq I_{\mu/ij}^0$.

\bigskip
Let $R=\Q[x_1,x_2,\ldots,x_n]$ be the polynomial ring in $n$ variables. The
inclusion above shows
that
  $$\dim\Big(R\big/ I^0_{\mu/ij}\Big)\  \le\  \dim\Big(R\big/
\widetilde{I}^0_{\mu/ij}\Big).
$$
We now present a reduction algorithm, modulo $\widetilde{I}_{\mu/ij}^0$,
that reduces any basis element
of $R$ as a linear combination of the basis elements of $R\big/
I^0_{\mu/ij}$. This will show that
   \begin{equation}\label{reduct}
   \dim\Big(R\big/ \widetilde{I}^0_{\mu/ij}\Big)\  \le\  \dim\Big(R\big/
{I}^0_{\mu/ij}\Big)
   \end{equation}
and conclude the proof of the Theorem~\ref{Imuij}.

{}From classical invariant theory, the set
  \begin{equation}\label{thebasis}
  \big\{ h_\lambda(\X) x_1^{\epsilon_1}x_2^{\epsilon_2}\cdots
   x_{n-1}^{\epsilon_{n-1}} \,:\,
    0\le\epsilon_i\le i-1 \big\}
   \end{equation}
forms a basis of $R$, as $\lambda$ runs through all partitions with parts
$\le n$. There are many proofs of
this fact. One that is more adapted to this context can be found in
\cite{nantel_adriano}, Theorem 3.2. We use the basis in \ref{thebasis} for our
reduction algorithm.

For $0\le \zeta\le\ell$, let $\nu(\zeta)$ be the partition of $n$ obtained
from $\mu$ by
removing the rightmost corner that is north-east of the cell
$(i+\zeta,j)$.

\centerline{
\epsffile{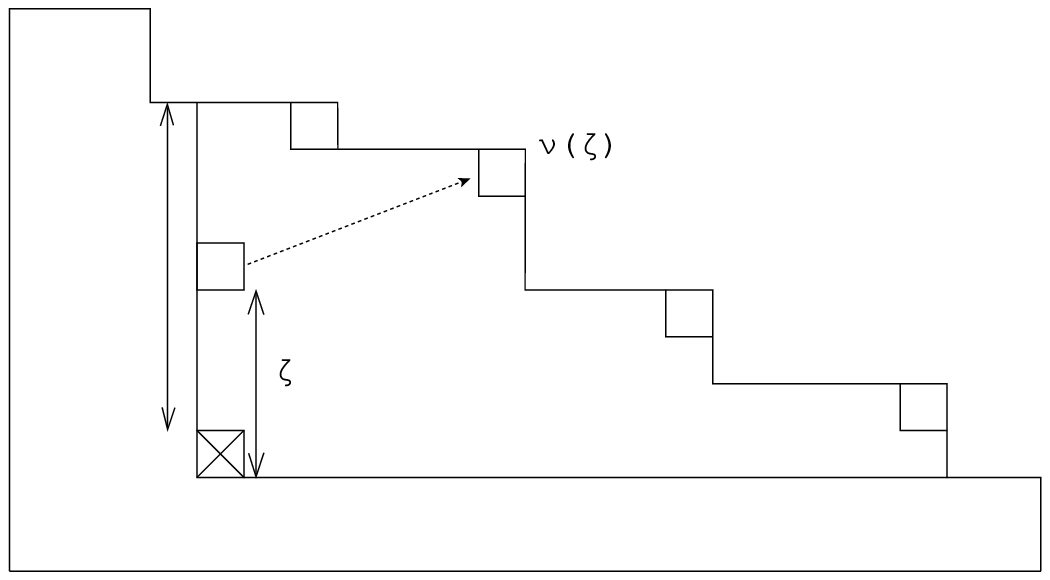}}

\centerline{\bf Figure 3}
\medskip

\noindent
Let ${\mathcal B}_{\nu(\zeta)}$ be class representative for a basis of
$R\big/ I^0_{\nu(\zeta)}$. Such a basis is given in
\cite{aval,nantel_adriano,GP}. Then we let
   \begin{equation}\label{basis}
   {\mathcal B}_{\mu/ij}=\  \sum_{\zeta=0}^\ell\  h_1^\zeta(\X) {\mathcal
B}_{\nu(\zeta)}
   \end{equation}
where $h_1^\zeta(\X) {\mathcal B}_{\nu(\zeta)}=\{h_1^\zeta(\X)P(\X) :P(\X)\in
{\mathcal B}_{\nu(\zeta)}\}$
and the sum indicates a (disjoint) union. We need the following lemmas.

\begin{lem} \label{star} For $0\le\zeta\le\ell$, we have
  $$ h_1^\zeta(\X)I^0_{\nu(\zeta)}\ \subseteq\  \widetilde{I}^0_{\mu/ij}\, + \,
\big\langle
     h_1^{\zeta +1}(\X)\big\rangle.
  $$
\end{lem}

\proof
Let $h_r(S)$ be a generator of $I^0_{\nu(\zeta)}$ where $|S|=k$ and
$r>\delta_k\big(\nu(\zeta)\big)-k$.
Let $(\alpha_\zeta,\beta_\zeta)$ be the coordinate of the corner cell such that
$\nu(\zeta)=\mu/\alpha_\zeta\beta_\zeta$. We note that
  $$
  \delta_k\big(\nu(\zeta)\big) =
   \begin{cases}
    \delta_k(\mu) & \hbox{ if $k\le\beta_\zeta $}\cr
    \delta_k(\mu)-1 & \hbox{ if $k>\beta_\zeta $}\cr
 \end{cases}
 $$
If $k\le\beta_\zeta$, then
$r>\delta_k\big(\nu(\zeta)\big)-k=\delta_k(\mu)-k$ and we have $h_r(S)\in
J_\mu^0\subseteq \widetilde{I}^0_{\mu/ij}$. Similarly, if $k>\beta_\zeta$ and
$r>\delta_k\big(\nu(\zeta)\big)-k+1=\delta_k(\mu)-k$, then again $h_r(S)\in
J_\mu^0\subseteq \widetilde{I}^0_{\mu/ij}$. We are left to assume that
$k>\beta_\zeta$ and
$r=\delta_k(\mu)-k$. For $\zeta=0$, we have
  $$I_{1^n}^0 =\langle h_r(\X) : r>0 \rangle\  \subseteq\
\widetilde{I}^0_{\mu/ij}\, + \, \big\langle
h_1(\X)\big\rangle,
$$
we can use Lemma~\ref{htoe} and
  $$h_r(S)\equiv (-1)^r e_r(\overline S) \mod \Big(\widetilde{I}^0_{\mu/ij}\,
+ \, \big\langle
h_1(\X)\big\rangle\Big),
$$
where $\overline S=\X-S$. In this case, $e_r(\overline S)\in
\widehat J_{\mu/ij}^0\subseteq \widetilde{I}^0_{\mu/ij}$ and thus
$h_r(S)\in\widetilde{I}^0_{\mu/ij}\, + \, \big\langle
h_1(\X)\big\rangle$. For $\zeta>0$, according to Lemma~\ref{htoe} again, we
have $h_r(S)\equiv (-1)^r
e_r(\overline S) \mod I_{1^n}^0$. In particular $h_1^\zeta(\X)h_r(S)\equiv
(-1)^r
h_1^\zeta(\X)e_r(\overline S) \mod h_1^\zeta(\X)I_{1^n}^0$. Since
$h_1^\zeta(\X)I_{1^n}^0\subseteq
\widetilde{I}^0_{\mu/ij}\, + \, \big\langle
     h_1^{\zeta +1}(\X)\big\rangle$, we have
  $$h_1^\zeta(\X)h_r(S)\equiv (-1)^r h_1^\zeta e_r(\overline S) \mod
\Big(\widetilde{I}^0_{\mu/ij}\, + \,
\big\langle h_1^{\zeta+1}(\X)\big\rangle\Big).
$$
In this case, $h_1(\X)e_r(\overline S)\in
\widehat J_{\mu/ij}^0\subseteq \widetilde{I}^0_{\mu/ij}$ and thus
$h_1^\zeta(\X)h_r(S)\in\widetilde{I}^0_{\mu/ij}\, + \,
\big\langle h_1^{\zeta+1}(\X)\big\rangle$.
\endproof

\begin{lem} \label{alg} Modulo $\widetilde{I}^0_{\mu/ij}$, any element of the
form  $h_\lambda(\X)
x_1^{\epsilon_1}x_2^{\epsilon_2}\cdots x_{n-1}^{\epsilon_{n-1}}$ with\break
    $0\le\epsilon_i\le i-1$ is a linear combination of elements in
${\mathcal B}_{\mu/ij}$.
\end{lem}

\proof We remark that $\big\langle
h_1^{\ell+1}(\X),h_2(\X),h_3(\X),\ldots\big\rangle\subseteq\widetilde{I}^0_{\mu/ij}$.
Hence
 \begin{equation}
  x_1^{\epsilon_1}x_2^{\epsilon_2}\cdots x_{n-1}^{\epsilon_{n-1}}\equiv 0
  \mod \widetilde{I}^0_{\mu/ij}
 \end{equation}
unless $h_\lambda(\X)=h_1^\zeta(\X)$ for $0\le \zeta\le\ell$. We then proceed
by induction on $\zeta$,
from $\zeta=\ell+1$ down to $\zeta=0$. The result is true for $\zeta>\ell$
since $h_1^{\ell+1}(\X)\equiv
0\mod
\widetilde{I}^0_{\mu/ij}$. For $\zeta\le\ell$ consider $h_1^\zeta(\X)
x_1^{\epsilon_1}x_2^{\epsilon_2}\cdots x_{n-1}^{\epsilon_{n-1}}$. We assume by induction that
modulo $\widetilde{I}^0_{\mu/ij}$, 
any $h_\mu x_1^{\eta_1}x_2^{\eta_2}\cdots x_{n-1}^{\eta_{n-1}}$ with 
 the number of parts equal to $1$ in $\mu$ is $\ge\zeta+1$ and   $0\le\eta_i\le i-1$, is a linear
combination of elements in
${\mathcal B}_{\mu/ij}$.
{}From our choice of ${\mathcal
B}_{\nu(\zeta)}$, there exists an element $A$ in the linear span of ${\mathcal
B}_{\nu(\zeta)}$ such that $x_1^{\epsilon_1}x_2^{\epsilon_2}\cdots
x_{n-1}^{\epsilon_{n-1}}\equiv A\mod
I^0_{\nu(\zeta)}$. Hence there is an element $B\in I^0_{\nu(\zeta)}$ such that
   \begin{equation}\label{recu}
   h_1^\zeta(\X) x_1^{\epsilon_1}x_2^{\epsilon_2}\cdots
x_{n-1}^{\epsilon_{n-1}} = h_1^\zeta(\X) A +
        h_1^\zeta(\X) B.
  \end{equation}
Here $h_1^\zeta(\X) A$ is in the linear span of $h_1^\zeta(\X) {\mathcal
B}_{\nu(\zeta)}\subseteq {\mathcal B}_{\mu/ij}$.  {}From Lemma~\ref{star},
$h_1^\zeta(\X) B\in
\widetilde{I}^0_{\mu/ij}\, + \, \big\langle  h_1^{\zeta +1}(\X)\big\rangle$. Hence
  $$ h_1^\zeta(\X) B\equiv h_1^{\zeta+1}C \mod \widetilde{I}^0_{\mu/ij},
 $$
where $C$ is an element of $R$. By our induction hypothesis
$h_1^{\zeta+1}C$ is in the linear span of
${\mathcal B}_{\mu/ij}$.
\endproof

We now recall an auxiliary result from~\cite{ABB}. For completeness
we will sketch the proof but the interested reader will find the complete details in the original
paper. We have that $\dim(R\big/I^0_{\mu/ij})=\dim(M_{\mu/ij}^0)$. Using this we have

\begin{prop} \label{dims}
  $$\dim(R\big/I^0_{\mu/ij})\ge \sum_{\zeta=0}^\ell \big| {\mathcal
B}_{\nu(\zeta)} \big| $$
\end{prop}

\proof[Sketch of proof] 
For this, we construct an explicit independent set of the right cardinality.
Recall that for $\nu$ a partition of $n$, a standard tableau $T$ of shape $\nu$ is a tableau that is
increasing both in rows and columns. We denote by $\RS{\nu}$ the set of standard tableaux of shape
$\nu$. We also associate to each entry
$j$, of a standard tableau
$T$,  a non-negative integer in the following manner. Let $(r_j,c_j)$ be the position of $j$ in $T$,
and let
$k$ be the largest entry of $T$, such that $c_k=c_j+1$ and $k<j$. We set $\gamma_T(j)=r_j-r_k.$
If there is no such $k$, set $\gamma_T(j)=r_j+1$. 

Recall that for $0\le \zeta\le\ell$ we denote by
$\nu(\zeta)=\mu/\alpha_\zeta\beta_\zeta$ the partition of
$n$ obtained from $\mu$ by removing $(\alpha_\zeta,\beta_\zeta)$, the  rightmost corner  that is
north-east of the cell
$(i+\zeta,j)$.  For $T$ a standard tableau, let
${\bf B}_T$ denote the set
   $${\bf B}_T=\big\{\ x_1^{m_1}x_2^{m_2}\cdots x_n^{m_n}\ |\  0\leq m_s\leq \gamma_T(s)\ \big\}.$$
A basis of $M^0_{\nu(\zeta)}$, is given by
   $${\mathcal B}_{\nu(\zeta)}=\big\{m(\partial \X)\Delta_T(\X)\ |\ T\in\RS{\nu(\zeta)},\   m(\X)\in
{\bf B}_T\ \big\},$$  
where $\Delta_T(\X)$ is defined in~\ref{garnirD}.

If $T$ is a standard tableau of shape $\nu(\zeta)$, and $0\leq u\leq \alpha_\zeta$ an
integer, we denote by
$T\!\uparrow_{u,\beta_\zeta}$ the tableau of shape $\mu/u\beta_\zeta$, such that
   $$T\!\uparrow_{u,\beta_\zeta}(r,c)=\left\{
 \begin{array}{cl}
   T(r,c)\hfill & \ \hbox{if} \ c\not=\beta_\zeta\ \hbox{or}\  r< u\,, \\
                 \\
   T(r-1,c)  & \ \hbox{if} \ c=\beta_\zeta\ \hbox{and}\ r> u \,.
 \end{array}
 \right.$$
In other words, the tableau $T\!\uparrow_{u,\beta_\zeta}$
is obtained from $T$ by ``sliding'' upward by $1$ the cells in column $\beta_\zeta$ that are on or
above row $u$. For $\mu/ij$, we set
   \begin{equation}\label{base muij}
       \tilde{\mathcal B}_{\mu/ij}=\bigcup_{\zeta=0}^\ell
           {\mathcal A}_{i+\alpha_\zeta-\zeta,\,\beta_\zeta},
   \end{equation}
where
   \begin{equation}\label{partie base}
      {\mathcal A}_{u,\,\beta_\zeta}=\big\{ m(\partial\X)
     \Delta_{T\uparrow_{u,\,\beta_\zeta}}(\X)\ |\ 
         T\in \RS{\nu(\zeta)},\   m(\X)\in{\bf B}_T  \big\}.
   \end{equation}
We prove that~\ref{base muij} is an independent set,
using a downward recursive argument.
Using $h_1(\partial \X)\Delta_{\mu/i,j}(\X,\Y)=\Delta_{\mu/i+1,j}(\X,\Y)$
in~\ref{introexpand} we obtain the following. For $T$  a standard tableau of shape
$\nu(\zeta)=\mu/\alpha_\zeta\beta_\zeta$ and
$0\leq u\leq \alpha_\zeta$, we
have\break $h_1(\partial\X)\Delta_{T\uparrow_{u,\beta_\zeta}}(\X)=\Delta_{T\uparrow_{u+1,v}}(\X)$ if
$u<\alpha_\zeta$, and $0$ if $u=\alpha_\zeta$. It follows from Definition~\ref{partie base} that
  $$
      h_1(\partial\X)\,{\mathcal A}_{u,\beta_\zeta}=\left\{
   \begin{array}{cl}
     {\mathcal A}_{u+1,\beta_\zeta}\hfill & \ \hbox{if} \ u<\alpha_\zeta\,, \\
                 \\
         \{0\}  & \ \hbox{if} \ u=\alpha_\zeta \,.
       \end{array}
      \right.
  $$
We deduce, from the linear independence of 
$B_{\nu(\zeta)}={\mathcal A}_{\alpha_\zeta,\beta_\zeta}$, that each ${\mathcal A}_{u,\beta_\zeta}$ is
independent. Applying
$h_1(\partial\X)$ in definition~\ref{base muij} we readily check that $\tilde{\mathcal
B}_{\mu/i+1,j}=h_1(\partial\X)\tilde{\mathcal B}_{\mu/i,j}$.
But we know that $h_1(\partial\X){\mathcal A}_{\alpha_\ell,\beta_\ell}=\{0\}$, and it is clear that
${\mathcal A}_{\alpha_\ell,\beta_\ell}$ is a subset of
$\tilde{\mathcal B}_{\mu/ij}$. By the induction hypothesis,
$\tilde{\mathcal B}_{\mu/i+1,j}$ is independent, and a counting argument forces the independence of
$\tilde{\mathcal B}_{\mu/ij}$.
\endproof

To conclude the proof of Theorem~\ref{Tmuij} we
use Lemma~\ref{alg} to show that
  $$ \dim(R\big/{\widetilde I}^0_{\mu/ij}) \le \big|{\mathcal B}_{\mu/ij} \big|
= \sum_{\zeta=0}^\ell \big|
{\mathcal B}_{\nu(\zeta)} \big|.$$ Then  the Equation~\ref{reduct} follows
from the Proposition~\ref{dims}.

\begin{cor} ${\mathcal B}_{\mu/ij}$ is a basis of $R\big/I^0_{\mu/ij}$.
\end{cor}

\section*{Acknowledgments}
We are indebted to Fran\c{c}ois Bergeron and Adriano Garsia for very
stimulating conversations and suggestions, and we thank Mary-Anne McLaughlin
for her careful reading. We also wish to thank the referee for his/her helpfull comments and
corrections.

\end{document}